\documentclass[a4paper,twoside,11pt]{article}

\usepackage[a4paper]{geometry}
\usepackage[latin1]{inputenc}
\usepackage[T1]{fontenc}

\usepackage{lipsum}
\usepackage{comment}
\usepackage{color}
\usepackage[english]{babel}
\usepackage{amssymb,amsmath,amsthm,amscd,amsbsy}
\usepackage{mathrsfs}
\usepackage{stmaryrd}
\usepackage{esint}
\usepackage{subfig}
\usepackage{tabularx} 
\usepackage{calc} 
\usepackage[pdftex]{graphicx} 
\usepackage[all]{xy}
\xyoption{v2}
\xyoption{2cell}
\UseAllTwocells

\usepackage{url}

\usepackage[pdfauthor={Christophe Prange},colorlinks=true, pdfstartview=FitV, linkcolor=blue, citecolor=blue, urlcolor=blue,pagebackref=false]{hyperref}

\usepackage{enumitem}

\usepackage{comment}

\usepackage{scalerel,stackengine}
\stackMath
\newcommand\reallywidehat[1]{%
\savestack{\tmpbox}{\stretchto{%
  \scaleto{%
    \scalerel*[\widthof{\ensuremath{#1}}]{\kern-.6pt\bigwedge\kern-.6pt}%
    {\rule[-\textheight/2]{1ex}{\textheight}}
  }{\textheight}%
}{0.5ex}}%
\stackon[1pt]{#1}{\tmpbox}%
}
\definecolor{darkgreen}{rgb}{0,0.5,0}
\definecolor{darkblue}{rgb}{0,0,0.7}
\definecolor{darkred}{rgb}{0.9,0.1,0.1}
\definecolor{lightblue}{rgb}{0,0.51,1}

\newtheorem{theorem}{Theorem}

\newcommand{\Z}{\mathbb{Z}}

\newcommand{\R}{\mathbb{R}}

\newcommand{\ep}{\varepsilon}

\renewcommand{\subset}{\subseteq}

\newcommand{\cu}{{\scaleobj{1.2}{\square}}}

\DeclareMathOperator{\dist}{dist}

\DeclareMathOperator{\supp}{supp}

\renewcommand{\bar}{\overline}

\begin{document}

\author{Christophe Prange\thanks{Universit\'e de Bordeaux, CNRS, UMR [5251], IMB, Bordeaux, France. \emph{e-mail:} \url{christophe.prange@math.u-bordeaux.fr}}}

\title{Infinite energy solutions to the Navier-Stokes equations in the half-space and applications}


\maketitle

\begin{abstract}
This short note serves as an introduction to the papers \cite{MMP17a,MMP17b}. These two works deal with the existence of mild solutions on the one hand and local energy weak solutions on the other hand to the Navier-Stokes equations in the half-space $\R^3_+$. We emphasize a concentration result for (sub)critical norms near a potential singularity. The contents of these notes were presented during the X-EDP seminar at IH\'ES in October 2017.
\end{abstract}

The aim of this short note is to provide an introduction to the results and techniques of two papers \cite{MMP17a,MMP17b} written in collaboration with Yasunori Maekawa (Kyoto University) and Hideyuki Miura (Tokyo Institute of Technology). Our work focuses on the Navier-Stokes equations
\begin{equation}
\label{e.nse}
\left\{ 
\begin{aligned}
& \partial_tu+u\cdot\nabla u-\Delta u+\nabla p = 0, \quad \nabla\cdot u=0  & \mbox{in} & \ (0,T)\times\R^3_+, \\
& u = 0  & \mbox{on} & \ (0,T)\times\partial \R^3_+. 
\end{aligned}
\right.
\end{equation}
in the half-space $\R^3_+$, with initial data $u_0$ locally uniformly Lebesgue integrable. Our paper about the linear theory \cite{MMP17a} treats the case of $\R^d_+$ in any dimension $d\geq 2$. Nevertheless, for the purpose of the present discussion we will restrict our attention to the case $d=3$. The literature on the mathematical study of the Navier-Stokes equations being so extensive, we will just mention some relevant works we know of related to our specific subject. For a more complete bibliography, the books of Lemari\'e-Rieusset \cite{lemariebook,nse21} are good references. We do not give any proof of the results stated in this note, but indicate some insights of the proofs or technical points. The proofs can be found in \cite{MMP17a,MMP17b}.

\smallskip

One of the most fundamental properties of the Navier-Stokes equations in $\R^3_+$ is the scaling: for all $\lambda>0$, if $(u,\nabla p)$ is a solution to \eqref{e.nse}, then $(u_\lambda,p_\lambda)$ defined by 
\begin{equation*}
u_\lambda(y,\tau):=\lambda u(\lambda y,\lambda^2\tau),\quad p_\lambda(y,\tau):=\lambda^2p(\lambda y,\lambda^2\tau),
\end{equation*}
for all $(y,\tau)\in\R^3_+\times (0,\lambda^{-2}T)$, is a solution to \eqref{e.nse} with initial data $u_{0,\lambda}:=\lambda u_0(\lambda\cdot)$. This property makes it possible to discriminate between functional spaces. Roughly speaking in subcritical spaces, the dynamics is driven by the linear part of the equation, in the critical spaces the linear and nonlinear parts are in balance, while in supercritical spaces the linear part is weak compared to the nonlinear part. In Section \ref{sec.loccon} we will study \eqref{e.nse} in subcritical and critical spaces, while the analysis in Section \ref{sec.lews} takes place in a supercritical space. Applications to blow-up criteria in terms of subcritical and critical norms will be provided.

\smallskip

Our work is focused on initial data $u_0\in L^q_{uloc,\sigma}(\R^3_+)$, for $q\in [1,\infty]$ i.e. (by definition)
\begin{align}
&u_0\in L^q_{loc}(\R^3_+),\quad \sup_{\eta\in\Z^3_+}\|u_0\|_{L^q(\eta+(0,1)^3)}<\infty,\label{e.condone}\\
\mbox{and}\quad &\int_{\R^3_+}u_0\cdot\nabla\varphi dx=0,\qquad\forall\varphi\in C^\infty_c(\overline{\R^3_+})\label{e.condtwo}.
\end{align}
The last condition implies that $\nabla\cdot u_0=0$ and $u_0\cdot e_3=0$ in the sense of distributions. Remark that $L^\infty_{uloc,\sigma}=L^\infty_\sigma$. For reasons related to the Navier-Stokes equations themselves (energy or scaling), we will need some local integrability on $u_0$, i.e. $q=2$ or $q\geq 3$ to define the solutions. We say that $f\in L^q_{uloc}(\R^3_+)$ for $q\in [1,\infty]$ if condition \eqref{e.condone} is satisfied. These local uniform Lebesgue spaces appear to have been introduced first by Kato \cite{K75} in the context of hyperbolic equations. There were later used in many different contexts such as linear parabolic equations \cite{ARCD04}, Ginzburg-Landau equations \cite{MS95,GV97}, reaction diffusion equations \cite{CD04}, water waves \cite{ABZ16}, boundary layers in fluids \cite{BGV08,DP14,DGV17} to cite just a few works. Their main interest lies in the fact that the functions have no quantitative decay at space infinity so that non trivial dynamics can arise from the equation itself, not due to exterior forcing. Moreover, they form a richer class than $L^\infty(\R^3_+)$ allowing for some singularities in the data. Indeed, the important scale critical function $u_0(x):=|x|^{-1}\in L^2_{uloc}(\R^3_+)$. The behavior can also be rather wild, as is shown by the following function in one space dimension: $u$ defined by $u(x):=n$, for $x\in(n,n+\frac1{n^2})$, zero elsewhere, belongs to $L^2_{uloc}(\R)$. Such functions become more and more concentrated like peaks at space infinity, while remaining uniformly locally in $L^q$. A variation on this example implies that $C^\infty_{b,\sigma}(\R^3_+)$ is not dense in $L^q_{uloc,\sigma}(\R^3_+)$ for $q\in[1,\infty]$. To rule out such behavior, we will sometimes consider the class
\begin{equation}\label{e.defmathcalL}
\mathcal L^q_{uloc,\sigma}(\R^3_+):=\overline{C^\infty_{c,\sigma}(\R^3_+)}^{L^q_{uloc}},\qquad q\in(1,\infty)
\end{equation}
We can characterize these functions (see \cite[Lemma A.4]{MMP17b}) as the functions $u_0$ of $L^q_{uloc,\sigma}(\R^3_+)$, which have some mild decay at infinity 
\begin{equation*}
\|u_0\|_{L^q(\eta+(0,1)^{2})}\rightarrow 0,\quad |\eta|\rightarrow\infty.
\end{equation*}
The notation $\overset{\circ}{E}_q$ is sometimes also used as in \cite{KS07,lemariebook}. Let us emphasize that the decay of the local norms at space infinity is not at all quantitative.

\smallskip

For initial data $u_0$ in the subcritical spaces $L^q_{uloc}(\R^3_+)$, $q\geq 3$, one can use the linear dynamics of the Stokes operator $\bf A$ to construct solutions to \eqref{e.nse} in a perturbative way. Treating the nonlinear term as a perturbation, we can write the following Duhamel formula
\begin{align}\label{e.mild.ns}
u(t)  = e^{-t{\bf A}} u_0 - \int_0^t e^{-(t-s){\bf A}} \mathbb{P} \nabla \cdot (u\otimes u ) d s  \,, \qquad t>0\,,
\end{align}
where $e^{-t{\bf A}}$ is the semigroup generated by the Stokes operator. Hence a fixed point of \eqref{e.mild.ns} classically yields a solution to the Navier-Stokes equations (see for instance \cite{Kato84} for pioneering work in the integrable case $L^q(\R^3)$, $q\geq 3$), which we call a mild solution. We will construct such solutions in Section \ref{sec.semi}.

\smallskip

For initial data in the supercritical space $L^2_{uloc}(\R^3_+)$, the smoothing of the linear Stokes operator is too weak to be able to solve \eqref{e.mild.ns}. Therefore, we rely on a priori bounds given by the (local) energy of the Navier-Stokes equations. If the initial data has finite energy i.e. $u_0\in L^2_{\sigma}(\R^3_+)$, the global energy of $u$ solving \eqref{e.nse} is bounded
\begin{equation}\label{e.globalen}
 \| u(t) \|_{L^2(\R^3_+)}^2 + 2\int_0^t \|\nabla u \|_{L^2 (\R^3_+)}^2 d s\leq \| u_0\|_{L^2(\R^3_+)}^2,
\end{equation}
for all $t\in(0,\infty)$. 
This was used to construct finite energy weak solutions, or so-called turbulent solutions, in the seminal works by Leray \cite{L34} and Hopf \cite{H51}. When the initial data is merely locally uniformly bounded, the energy of $u$ solving \eqref{e.nse} is not globally bounded (infinite energy), but only locally. We have the following local energy inequality
\begin{align}\label{e.locenineq} 
\begin{split}
& \| (\chi u)(t) \|_{L^2(\R^3_+)}^2 + 2\int_0^t \| \chi \nabla u \|_{L^2 (\R^3_+)}^2 d s \\
& \qquad \leq \| \chi u_0 \|_{L^2(\R^3_+)}^2+\int_0^t \langle |u|^2, \partial_s \chi^2 + \Delta \chi^2 \rangle_{L^2 (\R^3_+)}  + \langle u\cdot \nabla \chi^2, |u|^2 + 2 p \rangle_{L^2 (\R^3_+)} d s,
\end{split}
\end{align}
for all $\chi\in C^\infty_c([0,T)\times \overline{\R^3_+})$. The a priori bound \eqref{e.locenineq} at the difference of \eqref{e.globalen} involves the pressure in the flux of energy through $\partial\supp\chi$, which is a source of major complications. Such local energy weak solutions have been pioneered by Lemari\'e-Rieusset \cite{lemariebook} and are therefore sometimes named after him. We will investigate the existence of such solutions in the case of $\R^3_+$ in Section \ref{sec.lews}.

\smallskip

The main motivations of our work are twofold.

\noindent First, we aim at better understanding the role of the pressure as well as the effect of the boundary $\partial\R^3_+$ on the pressure. This is important for at least two reasons: the pressure has to be estimated in order to be able to use the local energy inequality \eqref{e.locenineq} to derive bounds on $u$, and the pressure cannot be eliminated using the Helmholtz-Leray projection on divergence-free fields which is unbounded on $L^\infty$ and a fortiori on $L^q_{uloc}$. Hence, we also develop an approach which circumvents the use of the Helmholtz-Leray projection.

\noindent Second, we want to use the existence theory for initial data barely locally uniformly Lebesgue integrable to investigate potential singularities of finite energy weak solutions to the Navier-Stokes equations. In particular, from the existence of mild solutions we derive immediately a concentration result for (sub)critical norms near the blow-up time, which gives a new direct proof of a recent result by Li, Ozawa and Wang \cite{LOW16}. Moreover, we apply the existence theory of local energy weak solutions to the proof of the blow-up of the scale critical norm $L^3(\R^3_+)$, thus recovering a result of Barker and Seregin \cite{BS15}.

\section{Two fundamental facts}
\label{sec.funda}

The first fact we want to stress is the existence of solutions which are driven by the pressure. For the non stationary Stokes system, we see that
\begin{equation*}
u(x,t):=f(t)\quad\mbox{and}\quad p(x,t):=-f'(t)\cdot x, 
\end{equation*}
is a solution in $\R^3$, while a solution in the half-space $\R^3_+$ is given by
\begin{equation}\label{e.parasit}
u(x,t):=(v_1(x_3,t),v_{2}(x_3,t), 0)\quad\mbox{and}\quad p(x,t):=-f(t)\cdot x', 
\end{equation}
where $f\in C^{\infty}_0((0,\infty);\R^{2})$ and $v(x_3,t)$ solves the heat equation $\partial_tv-\partial_d^2v=f$ with $v(0,t)=0$. Here and below we decompose the vector $x=(x',x_3)\in \R^3_+$ into the horizontal component $x'\in\R^2$ and the vertical component $x_3\in(0,\infty)$. We readily obtain that these special solutions are also solutions to the full Navier-Stokes system. Such solutions, for which the velocity is typically constant in space, are excluded in a finite energy setting. However they form a wide class of admissible solutions when one allows solutions which are non decaying at space infinity. Such solutions are dubbed ``parasitic solutions'' and have to be eliminated in order to get  an integral representation formula for the pressure in terms of the velocity (see Section \ref{sec.Liouville} for uniqueness theorems). 

\smallskip

The second fact is related to the presence of the boundary, namely to the computation of the pressure. In the whole space, the weak formulation of the Navier-Stokes equations immediately implies that
\begin{equation*}
-\Delta p=\nabla\cdot\nabla(u\otimes u)\qquad \mbox{in}\ \R^3.
\end{equation*}
If parasitic solutions are ruled out, one can thus get a representation formula for the pressure using the fundamental solution of the Laplacian. The pressure coincides with the Helmholtz pressure. In the half-space, the situation is more involved. Indeed the pressure solves the following Neumann problem
\begin{equation*}
\left\{ 
\begin{aligned}
& -\Delta p = \nabla\cdot(\nabla\cdot(u\otimes u))  & \mbox{in} & \ \R^3_+, \\
& \nabla p\cdot e_3=\gamma|_{x_3=0}\Delta u_3-\gamma|_{x_3=0}\nabla\cdot(u\otimes u)\cdot e_3 
& \mbox{in} & \ \partial\R^3_+.\\
\end{aligned}
\right.
\end{equation*}
The pressure is now a sum of the Helmholtz pressure and a pressure due to the boundary, which is called the ``harmonic pressure'' by several authors. This harmonic pressure $p_{harm}$ is harmonic in $\R^3_+$ and has Neumann boundary data $\gamma|_{x_3=0}\Delta u_3$. In the half-space $\R^3_+$, we obtain an explicit representation formula for the harmonic pressure via the formula
\begin{equation*}
p=P(x_3)(-\Delta)^{-\frac12}\gamma|_{x_3=0}\Delta u_3=-\frac{\nabla'\cdot}{(-\Delta)^{\frac12}}P(x_3)\gamma|_{x_3=0}\partial_{x_3=0}u',
\end{equation*}
where $P$ is the Poisson kernel for the half-space. 
Further formulas are given in Section \ref{sec.lews}. For more general domains $\Omega$ than the half-space, one usually lacks a formula for the harmonic pressure. In order to control the pressure, one relies on an a priori bound for $\dist(x,\partial\Omega)\nabla p$ in $L^\infty$. However such a bound is only available in certain domains called strictly admissible by Abe and Giga, such as bounded $C^3$ domains \cite{AG13}, exterior $C^3$ domains \cite{AG14}. In these two instances, the admissibility of the domain is proved using a blow-up argument. Further phenomena related to the harmonic pressure are investigated by Kang \cite{Kang05}.

\section{Resolvent estimates}
\label{sec.resol}

The first step in our analysis is the study of the linear evolution. In order to get estimates on the Stokes semigroup, we focus our attention on the stationary resolvent problem
\begin{equation}
\label{e.resol}
\left\{ 
\begin{aligned}
& \lambda v -\Delta v+ \nabla p = f, \quad  \nabla\cdot u=0  & \mbox{in} &\ \R^3_+, \\
& v = 0  & \mbox{on} & \ \partial\R^3_+. 
\end{aligned}
\right.
\end{equation}
for non localized and divergence-free data $f$. We are able to get estimates for \eqref{e.resol} for all $\lambda$ in the sector $S_{\pi-\ep}$ defined by
\begin{equation*}
S_{\pi-\ep}:=\{\rho e^{i \theta}:\ \rho>0,\, \theta\in[-\pi+\ep,\pi-\ep]\}\subset\mathbb C.
\end{equation*}
This in turn enables to get estimate in long time for the Stokes semigroup, see Section \ref{sec.semi} below. Our result is stated in the following theorem.

\begin{theorem}[Resolvent estimates, {\cite[Theorem 1]{MMP17a}}]
\label{theo.fresolvent}
Let $1 < q \leq \infty$, $\ep >0$. Let $\lambda$ be a complex number in the sector $S_{\pi-\ep}$.  Let $f\in L^q_{uloc,\sigma}(\R^3_+)$. Then there exist $C(\ep,q)<\infty$ (independent of $\lambda$) and  a solution $(v,\nabla p)\in L^q_{uloc}(\R^3_+) \times L^1_{uloc} (\R^3_+)$ to \eqref{e.resol} in the sense of distributions such that 
\begin{align} 
|\lambda| \| v\|_{L^q_{uloc}} + |\lambda|^\frac12 \| \nabla v \|_{L^q_{uloc}} & \leq C\| f\|_{L^q_{uloc}},\label{est.prop.fresolvent.1}\\
\| \nabla^2 v \|_{L^q_{uloc}} + \| \nabla p \|_{L^q_{uloc}} & \leq C (1+e^{-c|\lambda|^\frac12} \log |\lambda| )\| f\|_{L^q_{uloc}},\quad\mbox{for}\quad q\neq \infty.\label{est.prop.fresolvent.2}
\end{align}
Moreover, for $1 \leq q < p\leq \infty$ satisfying $\frac1q-\frac1p<\frac13$, there exists a constant $C(\ep,q,p)<\infty$ (independent of $\lambda$) such that
\begin{align}
\| v\|_{L^p_{uloc}} & \leq C | \lambda|^{-1}(1+|\lambda|^{\frac{3}{2}(\frac1q-\frac1p)}) \| f \|_{L^q_{uloc}}, \label{est.prop.fresolvent.3}\\
\| \nabla v\|_{L^p_{uloc}} & \leq C |\lambda|^{-\frac12}( 1+ |\lambda|^{\frac{3}{2}(\frac1q-\frac1p)}) \| f \|_{L^q_{uloc}}.\label{est.prop.fresolvent.4}
\end{align} 
\end{theorem}

Uniqueness of $v$ will be discussed in Section \ref{sec.Liouville}. There are issues specific to the locally uniform integrability framework, related to the existence of the parasitic solutions exhibited in Section \ref{sec.funda}.

\smallskip

Our main source of inspiration is the work of Desch, Hieber and Pr\"uss \cite{DHP01}, which is concerned with the Stokes resolvent problem in $L^\infty(\R^3_+)$. The techniques of the classical work of Farwig and Sohr \cite{FS94} for the Stokes resolvent problem in $L^q$, $q\in(1,\infty)$ based on the H\"ormander-Michlin theorem are neither applicable in the $L^\infty$ setting nor in the $L^q_{uloc}$ setting. The idea put forward in \cite{DHP01} is to decompose the solution $v$ to \eqref{e.resol} into a Dirichlet-Laplace part and a nonlocal part due to the pressure $v=v_{D.L.}+v_{nonloc}$.

\smallskip

The main advantage of the half-space $\R^3_+$ is the fact that one can get explicit formulas using the Fourier transform in the tangential variables. We denote by $v'$ the first two components of $v$ and $v_3$ the third component of $v$ so that $v=(v',v_3)$; in the same way we decompose $f=(f',f_3)$. We have for solutions which decay at space infinity, for all $\xi\in\R^2$, $y_3>0$,
\begin{align*}
\widehat{v_{D.L.}}(\xi,y_3)=\ &\frac{\xi}{2\omega_\lambda(\xi)} \int^{\infty}_0 (e^{-\omega_\lambda(\xi)|y_3-z_3|} -  e^{-\omega_\lambda(\xi)(y_3+z_3)})
\widehat{f}(\xi,z_3)dz_3\\
\widehat{v_{nonloc}'}(\xi,y_3)=\ &-\frac{i\xi}{2\omega_\lambda(\xi)} \int^{\infty}_0 (e^{-\omega_\lambda(\xi)|y_3-z_3|} -  e^{-\omega_\lambda(\xi)(y_3+z_3)})e^{-|\xi|z_3}\widehat{p_0}(\xi)dz_3\\
\widehat{v_{nonloc,3}}(\xi,y_3)=\ &\frac{|\xi|}{2\omega_\lambda(\xi)} \int^{\infty}_0 (e^{-\omega_\lambda(\xi)|y_3-z_3|} -  e^{-\omega_\lambda(\xi)(y_3+z_3)})e^{-|\xi|z_3}
\widehat{p_0}(\xi)dz_3,
\end{align*}
where $\omega_\lambda(\xi):=\sqrt{\lambda+|\xi|^2}$ and for $\xi\neq 0$
\begin{equation*}
\widehat{p_0}(\xi)=-\frac{\omega_\lambda(\xi)+|\xi|}{|\xi|}\int_0^\infty e^{-\omega(\xi)z_3}\widehat{f_3}(\xi,z_3)dz_3.
\end{equation*}
From these expressions it clearly appears that the nonlocal part is originating from the pressure. One further uses the fact that $f$ is divergence-free and that $f_3(\xi,0)=0$ to get a more tractable expression for the nonlocal part of the multiplier. We end up with
\begin{equation}\label{e.forhatvnl}
\widehat{v_{nonloc}}(\xi,y_3)\simeq \frac 1\lambda (e^{-|\xi|y_3}-e^{-\omega_\lambda(\xi)y_3})\frac{\xi\otimes\xi}{|\xi|}\int_0^\infty e^{-\omega_\lambda(\xi)z_3}\widehat{f'}(\xi,z_3)dz_3 
\end{equation}
where $\simeq$ means that the exact expression is a finite sum of similar terms. The term written in the right hand side of \eqref{e.forhatvnl} is the most singular among those one has to handle. 

\smallskip

Since there is no easy characterization of uniform local Lebesgue spaces $L^q_{uloc}$ on the Fourier side, we need to first derive pointwise estimates for the kernels in physical space associated to these multipliers. Then, we can proceed with the estimates for $v$ and its derivatives. For both tasks, the nonlocal part $v_{nonloc}$ deserves the most attention, the Dirichlet-Laplace part being more standard to handle. Hence we focus on the nonlocal part in the discussion below. 

\smallskip

Let us consider the kernel $s$ defined by
\begin{equation*}
s_\lambda (y',y_3,z_3) := \frac1\lambda \int_{\R^{d-1}} e^{iy'\cdot \xi} 
 \big ( e^{-|\xi| y_3} -  e^{-\omega_\lambda(\xi)y_3} \big )
\,e^{-\omega_\lambda(\xi)z_3} \frac{\xi\otimes\xi}{|\xi|}d\xi.
\end{equation*}
In view of estimating the action of $s_\lambda$ on non decaying functions, we need to be very careful both about the singularity of the kernel near $y=0$, and about the decay at large scales. By scaling arguments, one can focus on the case $|\lambda|=1$. From the factor $e^{-\omega_\lambda(\xi)z_3}$ we always gain exponential decay in the $z_3$ direction. Indeed, there exists $c(\ep)>0$ such that for all $\lambda\in S_{\pi-\ep}\cap\{|\cdot|=1\}$, for all $\xi\in\R^2$, for all $z_3>0$,
\begin{equation*}
|e^{-\omega_\lambda(\xi)z_3}|\leq e^{-cz_3}.
\end{equation*}
Thus we concentrate on the decay in $y'$ and $y_3$. We discuss between $y_3>1$ and $0<y_3\leq 1$. The first case is straightforward, because $y$ is away from the singularity at $0$. The second case, when $y$ is close to the boundary, requires to be cautious about the singularity near $y'=0$. We use more structure of the multiplier. For low frequencies $|\xi|\leq 1$, we use the fact that
\begin{equation*}
|e^{-|\xi|y_3}-e^{-\omega_\lambda(\xi)y_3}|\leq Cy_3,
\end{equation*}
and the corresponding bounds on the derivatives. For high frequencies $|\xi|>1$ we take advantage of extra smallness of the factor $e^{-|\xi| y_3} -  e^{-\omega_\lambda(\xi)y_3}$. Indeed, there exists $C>0$ such that for all $|\xi|>1$,
\begin{equation*}
\big|(e^{-|\xi|y_3}-e^{-\omega_\lambda(\xi)y_3})\tfrac{\xi\otimes\xi}{|\xi|}\big|\leq Cy_3e^{-|\xi|y_3},
\end{equation*}
together with parallel bounds for derivatives. In the end, we obtain (see \cite[Proposition 3.5]{MMP17a}) that there exist $C(\ep)>0$ such that for all $\lambda\in S_{\pi-\ep}$, $y'\in\R^2$, $z_3,\, y_3>0$,
\begin{equation}\label{e.pointwiseker}
|s_\lambda (y',y_3,z_3)|\leq \frac{C y_3}{(y_3 + z_3 +|y'|)^2}\frac{e^{-c |\lambda|^{\frac 12}z_3}}{\big (1+|\lambda|^{\frac 12}(y_3+z_3 + |y'|)\big ) \big (1+|\lambda|^{\frac 12}(y_3+z_3)\big ) }
\end{equation}
along with analogous bounds for derivatives. We remark that the pointwise bounds on the kernels derived in the paper \cite{MMP17a} did not previously appear in the work \cite{DHP01}. Notice also that the first factor in the right hand side of \eqref{e.pointwiseker} decays exactly as the fundamental solution of the Laplacian in $\R^3$. The additional decay achieved for $\lambda\in S_{\pi-\ep}$ is provided by the last two factors. This decay is crucial when it comes to estimating the action of $s_\lambda$ on $L^q_{uloc}$. 

\smallskip

Let us now turn to the estimates for $v$. For simplicity, we just sketch the estimate of 
\begin{equation*}
I[f']:=\int_{\R^2}\int_0^\infty s_\lambda(y'-z',y_3,z_3)f'(z',z_3)dz_3 dz'.
\end{equation*}
for $f'\in L^q_{uloc}(\R^3_+)$. We aim at estimating $\|I[f']\|_{L^p((0,1)^2\times (\eta_3,\eta_3+1))}$. We bound the integral using Young's inequality for convolutions in the tangential direction for $\tfrac1p=\tfrac1s+\tfrac1q-1$,  
\begin{equation*}
\|I[f'](\cdot,y_3)\|_{L^p((0,1)^2)}\leq I_1+I_2,
\end{equation*}
where the terms in the right hand side correspond respectively to the local part around the singularity at $0$
\begin{equation*}
I_1:=\sum_{\max |\alpha'_i|\leq 2,~\max |\alpha_i^{'}+ \beta_i^{'}|\leq 2}\int_0^\infty \|  r'_\lambda (\cdot, y_3,z_3) \|_{L^s(\alpha'+(0,1)^2)} \|  f' (\cdot, z_3) \|_{L^q(\beta'+(0,1)^2)} d z_3
\end{equation*}
and to the part away from the singularity
\begin{equation*}
I_2:=\sum_{\max |\alpha'_i|\geq 3,~\max |\alpha_i^{'}+ \beta_i^{'}|\leq 2}\int_0^\infty \|  r'_\lambda (\cdot, y_3,z_3) \|_{L^s(\alpha'+(0,1)^2)} \|  f' (\cdot, z_3) \|_{L^q(\beta'+(0,1)^2)} d z_3.
\end{equation*}
The estimate of $I_1$ follows from the observation that
\begin{align*}
I_1  \leq\ & C \int_0^{1} \| s_\lambda (\cdot, y_3, z_3) \|_{L^s(\R^{2})} \| f' (\cdot, z_3)\|_{L^q(\{|z'|\leq 8\})}  d z_3\\
&  \quad + C \sum_{n=1}^\infty \left ( \int_n^{n+1} \| s_\lambda (\cdot, y_3, z_3) \|_{L^s(\R^{2})}^{q'} d z_3 \right )^\frac{1}{q'} \| f' \|_{L^q_{uloc}}\\
=:\ & I_{1,1} + I_{1,2}
\end{align*}
and from the pointwise estimate \eqref{e.pointwiseker} 
\begin{equation*}
\|s_\lambda(\cdot,y_3,z_3)\|_{L^s(\R^{2})}\leq \frac{C \, e^{-c |\lambda|^\frac12 z_3}}{|\lambda|^\frac12 \big (1+|\lambda|^\frac12 (y_3+z_3)\big ) (y_3+z_3)^{2(\frac1q-\frac1p)}}.
\end{equation*}
When $p>q$ one gains extra decay in the vertical direction which enables to prove the bound
\begin{equation*}
\|I_1\|_{L^p((0,1)^2\times (\eta_3,\eta_3+1))}\leq C|\lambda|^{-1}(1+|\lambda|^{\frac32(\frac1q-\frac1p)})\|f'\|_{L^q_{uloc}},
\end{equation*}
even in the case when $q=1$. The case $p=q$ is slightly more involved as far as $I_{1,1}$ is concerned. We interpolate between the bound
\begin{equation*}
\left\|Tf'\right\|_{L^{1,\infty}_{y_3}(0,\infty)}\leq C|\lambda|^{-1}\|f'\|_{L^1_{z_3}(0,1)}
\end{equation*}
and the bound
\begin{equation*}
\left\|Tf'\right\|_{L^\infty_{y_3}(0,\infty)}\leq C|\lambda|^{-1}\|f'\|_{L^\infty_{z_3}(0,1)}
\end{equation*}
where $T$ is the operator defined by
\begin{equation*}
(T f' )(y_3)  = \int_0^1 \frac{e^{-c |\lambda|^\frac12 z_3}}{|\lambda|^\frac12(1+|\lambda|^\frac12 (y_3+z_3))} f' (z_3) d z_3.
\end{equation*}
Therefore, we get the bound 
\begin{equation*}
\|I_1\|_{L^q((0,1)^2\times (\eta_3,\eta_3+1))}\leq C|\lambda|^{-1}\|f'\|_{L^q_{uloc}},
\end{equation*}
only in the case $q\in(1,\infty]$. The estimate of $I_2$ uses a similar decomposition of the integral in $z_3$, but is simpler since the singularity of the kernel at $0$ is not seen. It can be estimated for all $1\leq q\leq p\leq\infty$.

\smallskip

We finish this discussion of Theorem \ref{theo.fresolvent} by three further remarks. 

\noindent The exponent $q=p=1$ is excluded seemingly for technical reasons. The only obstruction lies in the analysis of the nonlocal term in the vicinity of the singularity at $y=0$; see the term $I_{1,1}$ above. Such obstructions are often seen in the study of singular integral operators. Nevertheless, for all the other terms, in particular the Dirichlet-Laplace part, the exponents $1=q=p$ are allowed. Though our paper does not provide a further investigation of this point, we mention now two related works in the $L^1$ framework. In the paper \cite{DHP01}, it is proved that there exists $f\in L^1_\sigma(\R^3_+)$ such that $v_{nonloc}\notin L^1(\R^3_+)$. In the paper \cite{Koz98}, Kozono proves that there is a more fundamental obstruction in the study of the Navier-Stokes flow in $L^1$. In an exterior domain, strong solutions of the Navier-Stokes equations with values in $L^1$ only exist if the force exerted by the fluid on the boundary of the exterior domain is zero. That said, the fact that $q=1$ is forbidden here is really an effect from the boundary.

\noindent Our second observation concerns the estimate \eqref{est.prop.fresolvent.4} for the second order derivatives. We were not able to remove the $\log|\lambda|$ singularity for small $|\lambda|$. This results in slightly worse estimates on the longtime dynamics of the second order derivatives for the unsteady Stokes system. Moreover, due to this lack of uniform bounds in $|\lambda|$ on second order derivatives of \eqref{e.resol}, we are unable to apply the technique of \cite{FS94} to get the well-posedness for the Stokes system, i.e. when $\lambda=0$. 

\noindent We conclude this part by mentioning that resolvent estimates in $L^\infty$ for the Stokes system were obtained in more general strictly admissible domains by Abe, Giga and Hieber \cite{AGH15}. Their method is based on localization, use of $L^p$ resolvent estimates and interpolation of the $L^\infty$ norm between the $L^p$ norm of the function and its gradient. As a consequence, they obtain $L^\infty$ bounds for $v$ and $\nabla v$, but not for $\nabla^2v$ nor the gradient of the pressure.

\section{Liouville theorems for the linear Stokes system}
\label{sec.Liouville}

In the paper \cite{MMP17a}, we prove the following uniqueness theorem up to parasitic solutions for the resolvent system \eqref{e.resol}.

\begin{theorem}[Liouville for steady Stokes, {\cite[Theorem 4]{MMP17a}}]\label{theo.liou}
Let $\lambda\in S_{\pi-\ep}$. Any solution $(u,\nabla p)\in L^1_{uloc}(\R^3_+)\cap L^1_{uloc}(\R^3_+)$ with $p\in L^1_{loc}(\R^3_+)$ to the resolvent problem \eqref{e.resol} with $f=0$ is a parasitic solution.
\end{theorem}

In other words, $u=(a'(x_3),0)$ and $p=D\cdot x'+c$, where $D\in \mathbb C^{2}$ and $a'=a'(x_3)\in\mathbb C^2$ is smooth, bounded and has trace equal to zero on $x_d=0$. It is easily seen that $a'_i$ for all $i\in\{1,2\}$ has to solve the following problem
\begin{equation*}
\lambda a'-\partial_3^2a'=D,\qquad a'|_{x_3=0}=0,
\end{equation*}
which is the resolvent problem for the one-dimensional Dirichlet-Laplacian.

\smallskip

Moreover, Theorem \ref{theo.liou} easily gives uniqueness of the solution of the resolvent problem constructed in Theorem \ref{theo.fresolvent} under one of the additional conditions 
\begin{equation}\label{e.decaypressure}
\|\nabla'p\|_{L^1(|x'|\leq 1,R<x_3<R+1)}\rightarrow 0\quad\mbox{when}\quad R\rightarrow\infty
\end{equation}
of decay of the pressure in the vertical direction or 
\begin{equation}\label{e.decayu}
\|u\|_{L^1(|x'-y'|\leq 1,1<x_3<2)}\rightarrow 0\quad\mbox{when}\quad y'\rightarrow\infty.
\end{equation}
of decay of the velocity in the tangential direction. 

\smallskip

Theorem \ref{theo.liou} for the resolvent problem leads to the following Liouville theorem for solutions to the unsteady Stokes system
\begin{equation}\label{e.stokes}
  \left\{
\begin{aligned}
 \partial_t u - \Delta u + \nabla p & = f , \quad \nabla \cdot u = 0  \qquad \mbox{in}~ (0,T)\times \R^3_+\,, \\
u & = 0\quad \mbox{on}~ (0,T)\times \partial\R^3_+,  \qquad u|_{t=0}=u_0 \quad \mbox{in}~\partial \R^3_+,
\end{aligned}\right.
\end{equation}
Let $u_0\in L^1_{uloc,\sigma}(\R^3_+)$ and $f\in L^1_{loc}(\overline{\R^3_+}\times (0,T))$. We call $(u,\nabla p)$ a weak solution to \eqref{e.stokes} if $(u,\nabla p)$ 
satisfies the equations in the sense of distributions, $u\in L^\infty((0,T);L^1_{uloc,\sigma}(\R^{3}_+))$, $p,\ \nabla p\in L^1_{loc}(\R^3_+\times (0,T))$, for all $\delta\in(0,T)$
\begin{equation}\label{e.condpress}
\sup_{x\in\R^3_+}\int_\delta^T\|\nabla p(\cdot,t)\|_{L^1(B(x,1)\cap\R^3_+)}<\infty
\end{equation}
and $u$ is weakly continuous in time, i.e. for all $\varphi\in C^\infty_c(\overline{\R^3_+})$
\begin{equation}\label{e.weakcvt}
t\mapsto \int_{\R^3_+}u(x,t)\cdot\varphi(x)dx
\end{equation}
belongs to $C^0([0,T))$.

\begin{theorem}[Liouville for unsteady Stokes, {\cite[Theorem 5]{MMP17a}}]\label{theo.liounsteady}
Any weak solution $(u,\nabla p)$ to \eqref{e.stokes} in the sense above with $u_0=f=0$ is a parasitic solution in the sense of \eqref{e.parasit}.
\end{theorem}

A Liouville theorem for infinite energy solutions of the Stokes system in the half-space has been worked out by Jia, Seregin and Sverak \cite{JSS12,JSS13}. Their result holds for weak solutions in $L^\infty(\R^3_+\times(-\infty,0))$, so-called bounded ancient solutions. The proof in \cite{JSS12} is based on the use of the Fourier transform in the tangential variable, while the proof in \cite{JSS13} uses a duality argument. Our result holds under the weaker integrability assumption $u\in L^\infty((0,T);L^1_{uloc,\sigma}(\R^{3}_+))$. The proof of Theorem \ref{theo.liounsteady} relies on a duality argument similar to the one used for the proof of the steady theorem, Theorem \ref{theo.liou}. The time dependence makes the proof more intricate on a technical level. Indeed one has to regularize and take care of the integrability in space and time. Let us emphasize that a proof based on duality is in principle more versatile than the one based on Fourier analysis. However, it requires a precise knowledge of the decay of the solution to a linear unsteady Stokes problem with localized or fast decaying source term. In our analysis, this decay is obtained from the resolvent estimates of Theorem \ref{theo.fresolvent}, which uses the special structure of the half-space to apply Fourier analysis in the tangential direction. In the analysis of \cite{JSS13}, the authors rely on the bounds for the Stokes Green kernel derived by Solonnikov \cite{Sol03}.

\smallskip

Theorem \ref{theo.liounsteady} is required to recover a  representation formula for the pressure of the local energy weak solutions constructed in Section \ref{sec.lews}. The integrability condition on the gradient of the pressure \eqref{e.condpress} imposed on the weak solutions is convenient for our purposes. It enables to compensate for the lack of decay and integrability of the test functions due in particular to the harmonic pressure. Moreover, it is satisfied for the pressure given by the representation formula in Section \ref{sec.lews}. It is possible that the condition \eqref{e.condpress} can be removed, but we did not carry out further investigations of this point.

\section{Linear and bilinear estimates for the Stokes semigroup}
\label{sec.semi}

Our approach is to get estimates on the Stokes semigroup by using the bounds for the resolvent problem \eqref{e.resol}. As is classical, the main tool to go from the estimates for the stationary resolvent problem to the estimates for the time-dependent Stokes semigroup is Dunford's formula. In that perspective, our method is in the spirit of the work of Farwig and Sohr \cite{FS94} for $L^q$, $q\in(1,\infty)$, Desch, Hieber and Pr\"uss \cite{DHP01} for $L^q$, $q\in (1,\infty]$ and Abe, Giga and Hieber \cite{AGH15} for $L^\infty$. Other methods are the compactness method, new in this context, pioneered by Abe and Giga \cite{AG13}, and the more classical one based on direct estimates of the Green function carried out by Solonnikov \cite{Sol03} or by Maremonti and Starita \cite{MS03}. In this vein, let us mention efforts to estimate the Green function of the Stokes system in the half-space, which originate in the work of Ukai \cite{Ukai87} and were extended by Cannone, Planchon and Schonbek \cite{CPS00} and Danchin and Zhang \cite{DZ14}.

\smallskip

The results of Section \ref{sec.resol} and Section \ref{sec.Liouville} enable to define the Stokes operator $\mathbf A$ realized in $L^q_{uloc,\sigma}(\R^3_+)$ for any $q\in(1,\infty]$ by the formula $(\lambda+\mathbf A)^{-1}=R(\lambda)$ for all $\lambda\in S_{\pi-\ep}$. Here $R(\lambda)$ is the resolvent operator associating to $f\in L^q_{uloc,\sigma}(\R^3_+)$ the unique solution to \eqref{e.resol} satisfying the estimates of Theorem \ref{theo.fresolvent} and the uniqueness condition of Theorem \ref{theo.liou}. The domain $D(\mathbf A)$ can be described for finite $q$:
\begin{equation*}
D(\mathbf A)=\{u\in L^q_{uloc,\sigma}(\R^3_+):\ \nabla^\alpha u\in L^p_{uloc},\ \alpha=1,2,\ u=0\ \mbox{on}\ \partial\R^3_+\}.
\end{equation*}
This contrasts with the Stokes operator in $L^\infty_\sigma(\R^3_+)$, where we have no description of the domain of $\mathbf A$, because we are lacking an estimate such as \eqref{est.prop.fresolvent.2} in that case. The following result is one of our main theorems.

\begin{theorem}[Stokes semigroup, {\cite[Theorem 2]{MMP17a}}]\label{theo.semi}
For all $q\in(1,\infty]$, the Stokes operator $\mathbf A$ generates a bounded analytic semigroup $e^{-t\mathbf A}$ in $L^q_{uloc,\sigma}(\R^3_+)$.
\end{theorem}

Estimates for derivatives in time and space of $e^{-t\mathbf A}f$ for $f\in L^q_{uloc,\sigma}(\R^3_+)$ are given in \cite[Section 5]{MMP17a}. In the case of $L^\infty_\sigma(\R^3_+)$, our result is more precise than the one of Abe and Giga \cite{AG13}. In this paper, the authors obtain bounds on the semigroup by relying on a compactness argument. However, their bounds are only for short time because of the nature of the contradiction argument. Notice that an important point in our theorem is the fact that the fact that the semigroup is bounded, which describes the long time behavior of the flow. Theorem \ref{theo.semi} is a consequence of Theorem \ref{theo.fresolvent}. The long time behavior in particular is due to the fact that we can bound the solution $v$ to the resolvent problem for any $\lambda$ arbitrarily close to $0$ in the sector $S_{\pi-\ep}$ for some $\ep>0$. Notice that since $D(\mathbf A)$ is not dense in $L^q_{uloc,\sigma}(\R^3_+)$ the Stokes semigroup is not known to be strongly continuous in $L^q_{uloc,\sigma}(\R^3_+)$. 

\smallskip

In order to study the nonlinear Navier-Stokes equations, in particular the existence of mild solutions, we need to estimate the Oseen kernel $e^{-t\mathbf A}\mathbb P\nabla\cdot$, where $\mathbb P$ denotes the Helmholtz-Leray projection on divergence-free fields. The operator $\mathbb P$ being unbounded on $L^q_{uloc,\sigma}(\R^3_+)$ for any $q\in[1,\infty]$, we analyze the first-order operator $\mathbb P\nabla\cdot$ and combine it with the Stokes semigroup. Once again, bounds on the non stationary Oseen kernel are derived from bounds on $(\lambda+\mathbf A)^{-1}\mathbb P\nabla\cdot$ via the Dunford formula. We compute an explicit formula for $\mathbb P\nabla\cdot(u\otimes u)$. It is a sum of terms which correspond to derivatives in the physical space and of terms which are variations of
\begin{equation*}
\frac{\xi_\alpha\xi_\beta\xi_\gamma}{|\xi|}\int_0^{z_3}e^{-(z_3-s)|\xi|}\widehat{u_\alpha u_\beta}(\xi,s)ds
\end{equation*}
in Fourier space. The goal is then to understand the action of multipliers of the form $\frac{\xi_\alpha\xi_\beta\xi_\gamma}{|\xi|}e^{-t|\xi|}$ in physical space. It is classical when considering fluids with non localized data, see \cite{AKLN15,lemariebook,Serf95,TTY10} to cite a few works, to decompose between small scales and large scales in physical space. For instance in the case of the whole space $\frac{D_\alpha D_\beta D_\gamma}{-\Delta}$ is decomposed into
\begin{equation*}
\frac{D_\alpha D_\beta D_\gamma}{-\Delta}=D_\alpha D_\beta\big((1-\chi)\tfrac{D_\gamma}{-\Delta}\big)+D_\alpha D_\beta\big(\chi\tfrac{D_\gamma}{-\Delta}\big),
\end{equation*}
where $\chi$ is a cut-off function in physical space. Here we follow a similar insight, except that we cut-off in Fourier space instead of physical space. This enables on the one hand to keep the exponential decay in the vertical direction in high frequencies and on the other hand to consider fractional derivatives instead of the second-order derivatives $D_\alpha D_\beta$. Let $K=K(y',t)$ be the kernel associated to the multiplier $\frac{\xi_\alpha\xi_\beta\xi_\gamma}{|\xi|}e^{-t|\xi|}$ in Fourier space. For $\theta\in(0,1)$ and $\lambda\in S_{\pi-\ep}$, we decompose $K$ into
\begin{equation}\label{e.decK}
K=(-\Delta)^\frac{2-\theta}{2}K_{\theta,\geq|\lambda|^\frac12}+K_{\leq |\lambda|^\frac12},
\end{equation}
where the high frequency part is bounded as follows
\begin{equation}\label{e.hf}
|K_{\theta,\geq|\lambda|^\frac12}(y',t)|\leq \frac{Ce^{-t|\lambda|^\frac12}}{(|y'|+t)^{2+\theta}},\qquad\forall y'\in\R^2,\ t\in(0,\infty),
\end{equation}
and
\begin{equation}\label{e.lf}
|K_{\leq|\lambda|^\frac12}(y',t)|\leq \frac{C}{(|\lambda|^\frac12+|y'|+t)^{4}},\qquad\forall y'\in\R^2,\ t\in(0,\infty),
\end{equation}
for the low frequency part. Notice that in \eqref{e.hf} the singularity at $0$ is integrable when $\theta\in(0,1)$ and the decay in the vertical direction is fast, while for \eqref{e.hf} there is no singularity at $0$ and the kernel is integrable at infinity. The fractional derivative $(-\Delta)^\frac{2-\theta}{2}$ in \eqref{e.decK} is later combined to the operator $(\lambda+\mathbf A)^{-1}$ taking advantage of the resolvent estimates of Theorem \ref{theo.fresolvent}. All this results eventually in the following bilinear estimates, which underly the study of the perturbative, so-called mild solutions, to the Navier-Stokes system. 

\begin{theorem}[Bilinear estimates, {\cite[Theorem 3]{MMP17a}}]\label{theo.bili}  
Let  $1<q \leq p \leq \infty$ or $1\leq q<p\leq\infty$. Then for $\alpha=0,1$ and for all $t\in(0,\infty)$,
\begin{equation}\label{e.bilia} 
\| \nabla^\alpha  e^{-t {\bf A}} \mathbb{P} \nabla \cdot (u\otimes v) \|_{L^p_{uloc}}  \leq Ct^{-\frac{1+\alpha}{2}}  \big (t^{-\frac{3}{2}(\frac1q-\frac1p)} +1  \big ) \| u\otimes v \|_{L^q_{uloc}}, 
\end{equation}
and  
\begin{equation}\label{e.bilib} 
\| \nabla  e^{-t {\bf A}} \mathbb{P} \nabla \cdot (u\otimes v) \|_{L^q_{uloc}}  \leq Ct^{-\frac{1}{2}}  \big ( \| u\cdot\nabla v \|_{L^q_{uloc}} + \| v\cdot\nabla u \|_{L^q_{uloc}} \big ).
\end{equation}
\end{theorem}

Existence of mild solutions to the Navier-Stokes equations \eqref{e.nse} with initial data in the (sub)critical spaces $L^q_{uloc,\sigma}(\R^3_+)$ follows from Theorem \ref{theo.bili} in a classical way, see \cite[Proposition 7.1 and Proposition 7.2]{MMP17a}.

\section{Localized concentration of (sub)critical norms near blow-up time}
\label{sec.loccon}

It is known since the work of Leray \cite{L34} and the later work of Giga \cite{Giga86} that solutions to the Navier-Stokes equations in the whole space $\R^3$ blowing-up at time $T$ satisfy the following lower bound on the rate of blow-up
\begin{equation}\label{e.blowup}
\|u(\cdot,t)\|_{L^q(\R^3)}\geq C(T-t)^{-\frac12(1-\frac3q)}\quad\mbox{for}\quad q\in(3,\infty],\ t\in(0,T)
\end{equation}
with a constant $C(q)$. 
We recall that a solution is regular at the point $(x_0,t)$ for $x_0\in\R^3$ (resp. $\overline{\R^3_+}$) if there exists $r>0$ such that $u\in L^\infty(B(0,r)\times(t-r^2,t))$ (resp. $L^\infty((B(0,r)\cap \R^3_+)\times(t-r^2,t))$). If $(x_0,t)$ is not regular it is by definition singular. We say that $T$ is a blow-up time if it is the time of the first occurrence of a singular point. It took much longer to extend \eqref{e.blowup} to the critical case $q=3$. Escauriaza, Seregin and Sverak \cite{ESS03} proved that if $T$ is a blow-up time for $u$ then
\begin{equation*}
\limsup_{t\rightarrow T^-}\|u(\cdot,t)\|_{L^3(\R^3)}=\infty.
\end{equation*}
This result has been strengthened by Seregin \cite{Ser12} to $\|u(\cdot,t)\|_{L^3}\rightarrow\infty$ for $t\rightarrow T^-$ and extended to the case of the half-space by Barker and Seregin \cite{BS15}. We come back to the blow-up of scale critical norms in Section \ref{sec.lews} below. Notice that the contraposition of such blow-up criteria also gives conditions for regularity of the solutions. 

\smallskip

Recently, an interesting refinement of \eqref{e.blowup} has been worked out by Li, Ozawa and Wang \cite{LOW16}. They prove a localized version of the lower bound for the blow-up rate for $u$ associated to initial data $u_0\in L^\infty(\R^3)$ and blowing-up at time $T$. More precisely, for all $q\in[1,\infty]$, there exists constants $0<c,\, C<\infty$, there exists a sequence $t_k\rightarrow T^-$ and a sequence of points $x_k\in\R^3$ such that
\begin{equation}\label{e.locbluplq}
\|u(\cdot,t_k)\|_{L^q(|\cdot-x_k|\leq c\omega(t_k)^{-1})}\geq C\omega(t_k)^{1-\frac3q},
\end{equation}
where 
\begin{equation}\label{e.locbluplinfty}
\omega(t):=\|u(\cdot,t)\|_{L^\infty(\R^3)}\gtrsim (T-t)^{-\frac12}.
\end{equation}
Notice that the blow-up rate for the $L^\infty_\sigma$ norm holds for any $t\in(0,T)$. The proof is in two steps. The first consists in using the existence theory of mild solutions in $L^\infty_\sigma$ to obtain \eqref{e.locbluplinfty}. This is in the spirit of the technique of our paper, see below. The second step uses \eqref{e.locbluplinfty} and a frequency decomposition at a cut-off frequency roughly $\omega(t_k)$. This technique seems to originate from the work of Bourgain \cite{B99} on nonlinear Schr\"odinger equations. The authors show that 
\begin{equation}\label{e.projlf}
\|P_{\lesssim \omega(t_k)}u(\cdot,t_k)\|_{L^\infty}\geq \frac{\omega(t_k)}2
\end{equation}
for a well chosen sequence of times $t_k$. Here $P_{\lesssim \omega(t_k)}$ is the projector on low frequencies. The bound \eqref{e.locbluplq} follows from bounding from above the left hand side in \eqref{e.projlf} by $\|u(\cdot,t_k)\|_{L^q(|\cdot-x_k|\leq c\omega(t_k)^{-1})}$. In our work \cite{MMP17a} we realized that the existence of mild solutions in $L^q_{uloc,\sigma}(\R^3)$ (resp. $\R^3_+$) yields a direct proof of \eqref{e.locbluplq} in the case when $q\in[3,\infty]$.

\begin{theorem}[Localized minimal blow-up rate, {\cite[Corollary 1.1]{MMP17a}}]\label{theo.conc}
For all $q\geq 3$, there exists a positive constant $C(q)<\infty$ such that for all $0<T<\infty$, for all $u \in C((0,T);L^\infty_\sigma(\R^3_+))$ mild solution to \eqref{e.nse}, if $u$ blows up at $T$, then for all $t\in(0,T)$, there exists $x(t)\in\R^3_+$ with the following estimate
\begin{equation}\label{e.lowerbdblup}
\|u(t)\|_{L^q(|\cdot-x(t)|\leq \sqrt{T-t})} \ge \frac{C}{(T-t)^{\frac12(1-\frac 3q)}}.\\
\end{equation}
\end{theorem}

The proof of Theorem \ref{theo.conc} relies on the existence theory of mild solutions for $q\geq 3$. Therefore, it is restricted to $q\geq 3$, while \eqref{e.locbluplq} was proved for $q\in[1,\infty]$ in \cite{LOW16}. However, our result strengthens the result of Li, Ozawa and Wang when $q\neq\infty$ since the lower bound \eqref{e.lowerbdblup} holds for all $t\in(0,T)$ and not only along a subsequence. The key to the proof of \eqref{e.lowerbdblup} is the following rescaled existence result for mild solutions. In the spirit of \cite{MT06} which deals with $\R^3$, we prove, see \cite[Proposition 7.5]{MMP17a}, that for $\rho>0$ if the initial data $u_0$ is controlled
\begin{equation*}
\|u_0\|_{L^q_{uloc,(\rho)}}:=\sup_{\eta\in\Z^3_+}\|u_0\|_{L^q(\eta+(0,\rho)^3)}\leq C\rho^{\frac3q-1},
\end{equation*}
then there exists a unique solution $u\in L^\infty((0,T);L^q_{uloc,(\rho)})$ such that
\begin{equation*}
\sup_{0<t<T} \big ( \| u (\cdot,t) \|_{L^q_{uloc,(\rho)}} + t^\frac d{2q} \| u(\cdot,t) \|_{L^\infty} 
 \big ) \leq C\| u_0 \|_{L^q_{uloc,(\rho)}}
 \end{equation*}
with lifespan of at least $T\geq \rho^2$. This claim is a mere rescaling of the existence of mild solutions in $L^q_{uloc}=L^q_{uloc,(1)}$. Theorem \ref{theo.conc} follows from this fact by a simple contraposition argument. Our proof works in the case of the whole space as well as in the case of $\R^3_+$. The only tool we need is the Navier-Stokes scaling.

\smallskip

Let us finally emphasize that in the case in the case $q>3$, the lower bound \eqref{e.lowerbdblup} is stronger than the global one \eqref{e.blowup}. In the critical case $q=3$ on the contrary, the lower bound \eqref{e.lowerbdblup} is a much weaker result than the blow-up criteria of \cite{ESS03,Ser12,BS15}.

\section{Representation of the pressure and local energy weak solutions}
\label{sec.lews}

In this final section, we present the ideas required to construct weak solutions with initial data in the supercritical space $L^2_{uloc,\sigma}(\R^3_+)$. Our construction is in the vein of the work of Lemari\'e-Rieusset \cite{lemariebook} and Kikuchi and Seregin \cite{KS07} for the whole space $\R^3$. On the contrary to Leray-Hopf solutions for data in $L^2_\sigma$, the solutions we consider here are only locally in $L^2$ and so no global energy bound is expected. The only tool to control the solutions is the local energy inequality \eqref{e.locenineq}. A major difficulty of the local energy inequality is that it not only counts the energy that is dissipated, but also the energy which enters or leaves a unit cell. In particular, the energy flux involves the pressure. A first step is therefore to be able to compute the pressure via a representation formula. Apropos pressure, let us notice that the Leray-Hopf theory of finite energy weak solutions says nothing about the pressure. Extracting information about the pressure requires some non trivial work, especially in the case of domains with boundaries, see \cite{SvW86}. Returning to the topic of infinite (local) energy weak solutions, in the case of $\R^3$, the computation of the pressure was handled by various techniques. In \cite{KS07}, the authors include the representation formula in their definition of the local energy weak solutions. Hence, they have to check that the formula remains true for a weak limit of approximate solutions. In \cite{JS2}, the authors assume some mild, non quantitative, decay of the velocity at infinity, which enables to recover the formula for the pressure via a Liouville theorem.

\smallskip

In the case of $\R^3_+$, we made the decision to impose the following condition on the pressure of our local energy weak solutions
\begin{equation}\label{e.condpress}
\sup_{x\in \R^3_+} \left(\int_\delta^{T'} \| \nabla p \|_{L^\frac98 (B(x)\cap \R^3_+)}^\frac32 d t\right)^\frac23 <\infty\qquad\mbox{for all}\quad0<\delta<T'\leq T.
\end{equation}
If $T=\infty$, the condition is slightly modified $0<\delta<T'<\infty$. This condition is stable under weak limits and enables to apply our Liouville theorem, Theorem \ref{theo.liounsteady}, for the unsteady Stokes system. It could be that this condition on the pressure is not optimal and that the Liouville theorem could be generalized, but our aim was to have a self-contained and simple framework. Notice that the case of the half-space has a major difficulty as compared to the whole space, namely the fact that the pressure splits into a Helmholtz-Leray part and a non trivial harmonic part. Hence in $\R^3_+$ we have to assume the condition \eqref{e.condpress} in order to identify the pressure, while in $\R^3$ the mere mild decay of $u_0$ appears to be enough. Moreover, the expression for the pressure is far more involved in the case of $\R^3_+$ than $\R^3$.

\smallskip

In our paper \cite{MMP17b}, we subsequently define local energy weak solutions to \eqref{e.nse} for initial data in $\mathcal L^2_{uloc,\sigma}(\R^3_+)$. This mild decay automatically rules out parasitic solutions and ensures that one has a representation formula for the pressure. In the more general case when the initial data belongs to $L^2_{uloc,\sigma}(\R^3_+)$, we need an additional condition on the decay of the harmonic pressure at space infinity. Let $T\in(0,\infty]$. Roughly speaking, see \cite[Definition 1.1]{MMP17b} for a rigorous definition, a local energy weak solution $(u,\nabla p)$ is a weak solution to \eqref{e.nse}, such that $u\in L^\infty((0,T);\mathcal L^2_{uloc,\sigma}(\R^3_+))$ (or if $T=\infty$, $u\in L^\infty_{loc}([0,\infty);\mathcal L^2_{uloc,\sigma}(\R^3_+))$), $p\in L^\frac32_{loc}(\overline{\R^3_+}\times(0,T))$, such that the pressure satisfies \eqref{e.condpress}, such that the local energy inequality \eqref{e.locenineq} is satisfied for all $\chi\in C^\infty_c([0,T)\times \overline{\R^3_+})$, such that weak continuity in time \eqref{e.weakcvt} holds and such that $u(\cdot,t)$ strongly converges in $L^2(K)$ to $u_0$ for any compact set $K\Subset\overline{\R^3_+}$.

\smallskip

The representation formula for the pressure is obtained via the formulas from the linear theory presented in Section \ref{sec.resol} above and via Dunford's formula. The following principles serve as guidelines for the derivation of the representation formula. First we decompose the solution to \eqref{e.nse} into a linear part, which solves a Stokes system with initial data $u_0$, and a nonlinear part which solves a Stokes system with zero initial data and right hand side $-\nabla\cdot(u\otimes u)$. This gives a decomposition $p=p_{li}+p^{u\otimes u}$. For $p_{li}$ the Helmholtz-Leray part is of course zero and so $p_{li}$ is equal to its harmonic part. For $p^{u\otimes u}$ we have to deal with both the Helholtz-Leray part  and the harmonic part. Second, a general principle is to separate small scales from large scales as is done in the whole space $\R^3$. For large scales, we systematically take advantage of the fact that the pressure is defined only up to a time-dependent constant so that we can gain additional decay in the kernels. We do not enter into more details here because the derivation of the formulas is intricate and rigorously done in the paper \cite{MMP17b}. The formula for the pressure is derived in Section 2, while the identification of the pressure of the local energy weak solutions as defined above is done in Section 3 via the Liouville theorem, Theorem \ref{theo.liounsteady}. The following statement is the main result of the paper \cite{MMP17b}.

\begin{theorem}[Global local energy weak solutions, {\cite[Theorem 1]{MMP17b}}]\label{theo.globalweak} 
For any $u_0 \in \mathcal{L}^2_{uloc,\sigma} (\R^3_+)$ there exists a local energy weak solution $(u,\nabla p)$ to \eqref{e.nse} in $\R^3_+\times(0,\infty)$ with initial data $u_0$.
\end{theorem}

The solution constructed in Theorem \ref{theo.globalweak} is global in time. We are able to construct local in time local energy weak solutions for data in $L^2_{uloc,\sigma}(\R^3_+)$, that is with no decay at all at space infinity. As underlined above, in that case we need to add a condition on the decay of the pressure in the vein of \eqref{e.decaypressure} in order to recover the representation formula. This was done in \cite[Section 5]{MMP17b}. The construction of global in time solutions exploits two ideas. The first is to transfer the decay of the initial data to the solution itself at positive times. We prove, see \cite[Theorem 2]{MMP17b}, that for all $u_0\in\mathcal L^2_{uloc,\sigma}(\R^3_+)$, all local energy weak solution $u$ to \eqref{e.nse} on $\R^3_+\times(0,T)$ satisfies
\begin{equation}\label{e.decayatinfty}
\sup_{t\in (0,T)}\sup_{\eta\in\Z^3_+}\int_{\cu(\eta)}|\vartheta(\tfrac\cdot R) u(\cdot,t)|^2+\int_0^T\int_{\cu(\eta)}|\vartheta(\tfrac\cdot R)\nabla u|^2+\left(\int_0^T\int_{\cu(\eta)}|\vartheta(\tfrac\cdot R) u|^3\right)^\frac23\stackrel{R\rightarrow \infty}{\longrightarrow}0,
\end{equation}
for all $\delta\in(0,T)$, with $\vartheta$ a smooth cut-off such that $0\leq\vartheta\leq 1$, $\vartheta\equiv 0$ on $B(0,1)$ and $\vartheta\equiv 1$ on $\R^3\setminus B(0,2)$. It follows that for almost all positive time $t$, $u(\cdot,t)\in \mathcal L^4_{uloc,\sigma}(\R^3_+)$. The second idea is that one can decompose $u(\cdot,t)$ into a large part in $C^\infty_{c,\sigma}(\R^3_+)\subset L^2_{\sigma}$, for which we can rely on the existence theory of Leray-Hopf solutions, and a small part in $\mathcal L^4_{uloc,\sigma}(\R^3_+)$, for which we use the existence theory of mild solution. Such an idea which probably goes back to Calderon \cite{Cal90} (weak solutions in $L^q$ for $q\in(2,3)$) is now classical in the study of the Navier-Stokes system. 

\smallskip

Very recently, Maremonti and Shimizu uploaded a paper \cite{MS18} where they prove the existence of weak solutions in the half-plane $\R^2_+$ and the half-space $\R^3_+$ for non decaying data. Their solutions are global in time. For $d=2$ or $3$ and $q\in[d,\infty)$, the initial data $u_0$ belongs to $L^\infty(\R^d_+)\cap J^{1,q}_0(\R^d_+)$, where 
\begin{equation*}
J^{q}_0(\R^d_+)=\overline{C^\infty_{c,\sigma}(\R^d_+)}^{\dot{W}^{1,q}}
\end{equation*}
so that the gradient of $u_0$ belongs to $L^q$. A typical example of such data is a constant field plus a bump in $W^{1,q}_{0,\sigma}(\R^d_+)$. This is very different from our setting in Theorem \ref{theo.globalweak} for the existence of global in time solutions. Indeed we need some mild decay at infinity so that constant data are ruled out, but the initial data can be ``wild'' everywhere. The authors of \cite{MS18} prove in addition to the well-posedness a structure result for the solutions of \eqref{e.nse}. In the case $d=2$ and $q=2$, they decompose the solution as the solution of the linear Stokes equation with initial data $u_0$ and the solution to a perturbed nonlinear Navier-Stokes equation with zero initial data, for which they can apply the $L^2$ theory. For higher $q$ intermediate perturbed linear problems are involved in the decomposition, until one reaches the problem for which $L^2$ theory can be used; see \cite{BA18} for a similar construction.

\smallskip

We finish this note by commenting on blow-up criteria for scale critical norms. As we mentioned previously, Seregin proved in \cite{Ser12} that for a Leray-Hopf solution $u$ to the Navier-Stokes system in the whole space which blows-up at time $T$, we have
\begin{equation*}
\|u(\cdot,t)\|_{L^3(\R^3)}\rightarrow\infty\quad\mbox{when}\quad t\rightarrow T^-.
\end{equation*}
This result was later extended to Leray-Hopf solutions in $\R^3_+$ by Barker and Seregin \cite{BS15} through different techniques. We are able to recover the result of \cite{BS15} for the half-space by using the technique of \cite{Ser12}. So let us outline the method in \cite{Ser12} and point out where the theory of local energy weak solutions is useful. Assume by contraposition that $\|u(\cdot,t_k)\|_{L^3}$ remains bounded along a subsequence of times $t_k\rightarrow T^-$. Without loss of generality, assume that $(0,T)\in \overline{\R^3_+}\times (0,\infty)$ is a blow-up point. The idea is to consider the sequence $u^{(k)}$ defined by blowing-up around the singularity at $0$ 
\begin{equation*}
u^{(k)}(y,s):=\lambda_ku(\lambda_k y,\lambda_k^2 s+T)\quad\mbox{with}\quad \lambda_k:=\sqrt{\frac{T-t_k}{S}}
\end{equation*}
for a parameter $S$ taken later on sufficiently small. The point is then to study the convergence of $u^{(k)}$. Since the energy in \eqref{e.globalen} is supercritical with respect to the Navier-Stokes scaling, the energy of the $u^{(k)}$'s is not uniformly bounded in $k$. The limit $\bar u$ is a local energy weak solution. Two things have to be proved about $\bar u$. First $\bar u(\cdot,0)=0$. This is a consequence of the fact that $u(\cdot,T)\in L^3$. Second, since the $L^3$ norm is scaling invariant, $\|u^{(k)}(\cdot,-S)\|_{L^3}=\|u(\cdot,t_k)\|_{L^3}$ is uniformly bounded. This implies some decay of the local $L^2$ norms of $u^{(k)}$ at space infinity. The key is to transfer this decay to $u^{(k)}(\cdot,s)$ at almost all times $s\in(-S,0)$. The tool for this is an estimate in the vein of \eqref{e.decayatinfty}. Our main contribution for the half-space is in this estimate. It appears that an estimate such as the decay bound \eqref{e.decayatinfty} was missing to Barker and Seregin in the half-space. Hence they found another technique to circumvent this point. The last part of the proof is similar in all the works \cite{Ser12,BS15,MMP17b,BA18} originating from \cite{ESS03}. It consists in proving a Liouville type theorem for $\bar u$ using backward uniqueness and unique continuation for parabolic-type equations. For the details in $\R^3_+$, we refer to \cite[Section 7]{MMP17b}.

\section*{Acknowledgement}

The author would like to thank Yasunori Maekawa and Hideyuki Miura for reading an earlier version of this manuscript. The author acknowledges financial support from the French Agence Nationale de la Recherche under grant ANR-16-CE40-0027-01, as well as from the IDEX of the University of Bordeaux for the BOLIDE project.

\small
\bibliographystyle{abbrv}
\bibliography{xedp}

\end{document}